\documentclass[graybox]{svmult}

\usepackage{type1cm}        
\usepackage[bottom]{footmisc}

\usepackage{amsmath,amsfonts,amssymb}
\usepackage{enumitem}
\usepackage{tikz}
\usepackage{array}
\usepackage{makecell}
\newcolumntype{x}[1]{>{\centering\arraybackslash}p{#1}}
\usepackage{tikz}
\usepackage{todonotes}

\usepackage{algorithm}
\usepackage{algorithmic}

\newcommand{\R}{\mathbb{R}}
\newcommand{\N}{\mathbb{N}}
\newcommand{\yw}{\widetilde{y}}
\newcommand{\pw}{\widetilde{p}}

\newcommand{\by}{{\bf y}}
\newcommand{\bd}{{\bf d}}

\makeindex             

\begin{document}
	
	\title*{Nonlinear optimized Schwarz preconditioner for elliptic optimal control problems}
	\titlerunning{Nonlinear OSM preconditioning for optimal control}
	\author{Gabriele Ciaramella, Felix Kwok and Georg M\"uller}
	\authorrunning{G. Ciaramella, F. Kwok and G. M\"uller}
	\institute{G. Ciaramella \at Politecnico di Milano \email{gabriele.ciaramella@polimi.it}	
	\and 
	F. Kwok \at Universit\'e Laval \email{felix.kwok@mat.ulaval.ca}
	\and 
	G. M\"uller \at Universit\"at Konstanz \email{georg.mueller@uni-konstanz.de}}
	%
	%
	\maketitle
	
	\abstract*{We introduce a domain decomposition-based nonlinear preconditioned
iteration for solving nonlinear, nonsmooth elliptic optimal control
problems, with a nonlinear reaction term, $L^1$ regularization and box
constraints on the control function. The method is obtained by
applying semismooth Newton to the fixed-point equation of the parallel
optimized Schwarz iteration. As a proof of concept, numerical
experiments are performed on two subdomains, as well as on a
multi-subdomain test case. The results show that it is possible
to obtain substantial improvements in robustness and efficiency with the
new method, relative to semismooth Newton applied directly to the
full optimization problem, provided appropriate Robin parameters and
a good continuation strategy are chosen.}

\section{Introduction}
\vspace*{-4mm}
	
Consider the nonlinear optimal control problem
\begin{equation}\label{eq:OCP_nonl}
\begin{split}
\min_{y,u} &\; J(y,u):=  \frac{1}{2} \| y-y_d \|_{L^2}^2 + \frac{\nu}{2} \| u \|_{L^2}^2
+ \beta \| u \|_{L^1}, \\
\text{s.t. } &-\Delta y + c y + b \varphi(y) = f + u \text{ in $\Omega$}, \: y =0 \text{ on $\partial \Omega$}, \\
&u \in U_{\rm ad} := \{ v \in L^2(\Omega) \, : \, |v| \leq \bar u \text{ in $\Omega$} \},
\end{split}
\end{equation}
where $\|\cdot\|_{L^r}$ denotes the usual norm for $L^r(\Omega)$ with $1\leq r \leq \infty$,
the functions $y_d,f\in L^2(\Omega)$ are given, and the scalar parameters $b,c,\beta \geq 0$ and $\nu,\beta \geq 0$ are known.
Our model includes problems such as the simplified Ginzburg-Landau
superconductivity equation as well as inverse problems where $L^1$-regularization is used
to enhance sparsity of the control function $u$.
For simplicity, the domain $\Omega \subset \R^2$ is assumed 
to be a rectangle $(0,\widetilde{L})\times(0,\widehat{L})$.
The function $\varphi : \R \rightarrow \R$ is assumed to be of class $C^2$, with locally bounded
and locally Lipschitz second derivative and such that $\partial_y\varphi(y)\geq 0$.
These assumptions guarantee that the Nemytskii operator $y(\cdot)\mapsto \varphi(y(\cdot))$ is 
twice continuously Fr\'echet differentiable in $L^\infty(\Omega)$.
In this setting, the optimal control problem \eqref{eq:OCP_nonl}
is well posed in the sense that there exists a minimizer
$(y,u) \in X \times L^2(\Omega)$, with $X:=H^1_0(\Omega) \cap L^\infty(\Omega)$, cf.~\cite{troltzsch2010optimal,Herzog}. Our goal is to derive efficient nonlinear preconditioners 
for solving \eqref{eq:OCP_nonl} using domain decomposition techniques.

Let $(y,u) \in X \times L^2(\Omega)$ be a solution to \eqref{eq:OCP_nonl}.
Then there exists an adjoint variable $p \in X$
such that $(y,u,p)$ satisfies the system \cite[Theorem 2.3]{Stadler2009}

\begin{align*}
-\Delta y + c y + b \varphi(y) &= f + u &&\text{in $\Omega$ with $y =0$ on $\partial \Omega$}, \\
-\Delta p + c p + b \varphi'(y)p &= y-y_d &&\text{in $\Omega$ with $p =0$ on $\partial \Omega$}, \\
u &= \mu(p),
\end{align*}
where $\mu : L^\infty(\Omega) \rightarrow L^2(\Omega)$ is
\begin{equation}\label{eq:mu}
\begin{split}
\mu(p) = &\max(0,(- \beta - p)/\nu)+\min(0,(\beta - p)/\nu) \\
& -\max(0,-\bar{u}+(- p - \beta)/\nu) - \min(0,\bar{u}+(- p + \beta)/\nu).
\end{split}
\end{equation}
We remark that for $\beta=0$, the previous formula becomes
$\mu(p) = \mathbb{P}_{U_{\rm ad}}(- p / \nu)$,
which is
the usual projection formula that leads to the optimality condition $u=\mathbb{P}_{U_{\rm ad}}(- p / \nu)$;
see \cite{troltzsch2010optimal}.
Moreover, if $\beta=0$ with $\bar u =\infty$, one obtains that $\mu(p)=- p/\nu$,
which implies the usual optimality condition $\nu u + p=0$, 
where $\nu u + p$ is the gradient of the reduced
cost functional $\widehat{J}(u)=J(y(u),u)$ \cite{troltzsch2010optimal}.

Eliminating the control using $\mu(p)$, the first-order optimality system becomes
\begin{equation}\label{eq:system}
\begin{split}
-\Delta y + c y + b \varphi(y)   &= f + \mu(p) \text{ in $\Omega$ with $y =0$ on $\partial \Omega$}, \\
-\Delta p + c p + b \varphi'(y)(p) &= y-y_d    \quad \text{ in $\Omega$ with $p =0$ on $\partial \Omega$}.
\end{split}
\end{equation}
This nonlinear and nonsmooth system admits a
solution $(y,p) \in X^2$ \cite{Herzog,troltzsch2010optimal}.

\section{Optimized Schwarz method and preconditioner}
\vspace*{-4mm}

In this section, we introduce an optimized Schwarz method (OSM) for solving 
the optimality system \eqref{eq:system}.
\begin{figure}
\centering
    \begin{tikzpicture}[xscale=3.6,yscale=0.52]
%

      \begin{scope}[shift={(0,-0.3)}]
      \draw (0,0) -- (.6,0)-- (.6,2) -- (0,2) -- cycle node[] at (.3,1) {$\Omega_1$};
		\node[] at (.06,1) {$\Gamma_0$};      
      \node[] at (.53,1) {$\Gamma_1$};
      \node at (0.85,1) {\small{$\cdots$}};
      \draw [dashed] (.6,0) -- (1.1,0);
      \draw [dashed] (.6,2) -- (1.1,2);
      \draw (1.1,0) -- (1.7,0)-- (1.7,2) -- (1.1,2) -- cycle node[] at (1.4,1) {$\Omega_j$};
      \node[] at (1.21,1) {$\Gamma_{j-1}$};
      \node[] at (1.63,1) {$\Gamma_j$};
      \node at (1.95,1) {\small{$\cdots$}};
      \draw [dashed] (1.7,0) -- (2.2,0);
      \draw [dashed] (1.7,2) -- (2.2,2);
      \draw (2.2,0) -- (2.8,0)-- (2.8,2) -- (2.2,2) -- cycle node[] at (2.5,1) {$\Omega_N$};
      \node[] at (2.31,1) {$\Gamma_{N-1}$};
      \node[] at (2.73,1) {$\Gamma_{N}$};
      \draw [<->] (-.1,0) --(-.1,2) node [midway, left] {$\widehat{L}$};
      \draw [<->] (0.0,2.2) --(0.6,2.2) node [midway, above] {$L$};
      \draw [<->] (1.1,2.2) --(1.7,2.2) node [midway, above] {$L$};
      \draw [<->] (2.2,2.2) --(2.8,2.2) node [midway, above] {$L$};
      \end{scope}
      \end{tikzpicture}
    \caption{Non-overlapping domain decomposition.}
    \label{fig:GRID2DDN}
\end{figure}
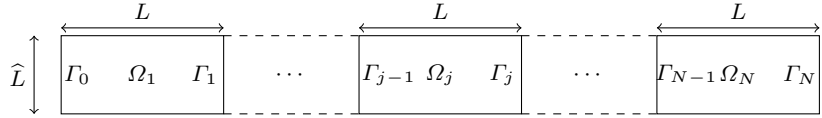
We consider the non-overlapping decomposition of $\Omega$ shown in Fig. \ref{fig:GRID2DDN}
and given by disjoint subdomains $\Omega_j$, $j=1,\dots,N$ such that $\overline{\Omega} = \cup_{j=1}^N \overline{\Omega}_j$.
The sets $\Gamma_j := \overline{\Omega}_j \cap \overline{\Omega}_{j+1}$,
$j=1,\dots,N-1$ are the interfaces. Moreover, we define $\Gamma_j^{\rm ext} := \partial \Omega_j \cap \partial \Omega$, $j=1,\dots,N$, which represent the external boundaries of the subdomains.
The optimality system \eqref{eq:system} can be written as a coupled system of 
$N$ subproblems defined on the subdomains $\Omega_j$, $j=1,\dots,N$, of the form
\begin{subequations}\label{eq:sub1}
\begin{align}
-\Delta y_j + c y_j + b \varphi(y_1)   &= f_j + \mu(p_j) &&\text{ in $\Omega_j$}, \label{eq:sub1_1} \\
-\Delta p_j + c p_j + b \varphi'(y_j)(p_j) &= y_j-y_{d,j}    &&\text{ in $\Omega_j$} \label{eq:sub1_2} \\
y_j &=0, \; p_j =0 &&\text{ on $\Gamma_j^{\rm ext}$}, \label{eq:sub1_3} \\
q \, y_j + \partial_{x} y_j &= q \,y_{j+1} + \partial_{x} y_{j+1} &&\text{ on $\Gamma_j$}, \label{eq:sub1_4} \\
 q \, p_j + \partial_{x} p_j &= q \, p_{j+1} + \partial_{x} p_{j+1} &&\text{ on $\Gamma_j$},  \label{eq:sub1_5} \\
 q \, y_j -\partial_{x} y_j &= q \,y_{j-1} -\partial_{x} y_{j-1} &&\text{ on $\Gamma_{j-1}$}, \label{eq:sub1_6} \\
q \, p_j -\partial_{x} p_j &= q \, p_{j-1} -\partial_{x} p_{j-1} &&\text{ on $\Gamma_{j-1}$},  \label{eq:sub1_7}
\end{align}
\end{subequations}
for $j=1,\dots,N$, where for $j\in\{1,N\}$ the boundary
conditions at $\Gamma_0$ and $\Gamma_N$, respectively, must be replaced with homogeneous Dirichlet conditions.
Here, $q>0$ is a parameter that can be optimized to improve
the convergence of the OSM; see, e.g, \cite{Gander2006,Chaouqui2018}.
The system \eqref{eq:sub1} leads to the OSM, which, for
a given 
$(y_j^{0},p_j^{0})_{j=1}^N$, is defined by solving the 
subdomain problems
\begin{subequations}\label{eq:OSM_sub1}
\begin{align}
-\Delta y_j^k + c y_j^k + b \varphi(y_1^k)   &= f_j + \mu(p_j^k) &&\text{ in $\Omega_j$}, \label{eq:OSM_sub1_1} \\
-\Delta p_j^k + c p_j^k + b \varphi'(y_j^k)(p_j^k) &= y_j^k-y_{d,j}    &&\text{ in $\Omega_j$} \label{eq:OSM_sub1_2} \\
y_j^k &=0, \; p_j^k =0 &&\text{ on $\Gamma_j^{\rm ext}$}, \label{eq:OSM_sub1_3} \\
q \, y_j^k + \partial_{x} y_j^k &= q \,y_{j+1}^{k-1} + \partial_{x} y_{j+1}^{k-1} &&\text{ on $\Gamma_j$}, \label{eq:OSM_sub1_4} \\
 q \, p_j^k + \partial_{x} p_j^k &= q \, p_{j+1}^{k-1} + \partial_{x} p_{j+1}^{k-1} &&\text{ on $\Gamma_j$},  \label{eq:OSM_sub1_5} \\
 q \, y_j^k -\partial_{x} y_j^k &= q \,y_{j-1}^{k-1} -\partial_{x} y_{j-1}^{k-1} &&\text{ on $\Gamma_{j-1}$}, \label{eq:OSM_sub1_6} \\
q \, p_j^k -\partial_{x} p_j^k &= q \, p_{j-1}^{k-1} -\partial_{x} p_{j-1}^{k-1} &&\text{ on $\Gamma_{j-1}$}  \label{eq:OSM_sub1_7}
\end{align}
\end{subequations}
for $k \in \N^+$.
Now, we use the OSM to introduce a nonlinear preconditioner by
setting $\by_j:=(y_j,p_j)$, $j=1,\dots,N$,
and defining the solution maps $S_j$ as
\begin{align*}
S_1(\by_2) &= \by_1 &&\text{solution to \eqref{eq:sub1} with $j=1$ and $\by_2$ given}, \\
S_j(\by_{j-1},\by_{j+1}) &= \by_j &&\text{solution to \eqref{eq:sub1} with $2\leq j \leq N-1$ and $\by_{j-1},\by_{j+1}$
given},\\
S_N(\by_{N-1}) &= \by_N &&\text{solution to \eqref{eq:sub1} with $j=N$ and $\by_{N-1}$ given}.
\end{align*}
Hence, using the variable $\by=(\by_1,\dots,\by_N)$,
we can rewrite \eqref{eq:sub1} as
\begin{equation}\label{eq:root_prec}
\mathcal{F}_{\rm P}(\by) = 0\text{, } \quad \text{where} \quad
\mathcal{F}_{\rm P}(\by) :=
\begin{small}
\begin{bmatrix}
\by_1 - S_1(\by_2) \\
\by_2 - S_2(\by_1,\by_3) \\
\vdots \\
\by_{N-1} - S_{N-1}(\by_{N-2},\by_N) \\
\by_N - S_N(\by_{N-1}) \\
\end{bmatrix} 
\end{small}.
\end{equation}
This is the nonlinearly preconditioned form of \eqref{eq:system}
induced by the OSM \eqref{eq:sub1}-\eqref{eq:OSM_sub1}, to which
we can apply a generalized Newton method.
For a given initialization $\by^0$,
a Newton method generates a sequence $(\by^k)_{k\in\N}$
defined by
\begin{equation}\label{eq:SSN_prec}
\text{solve } \quad D\mathcal{F}_{\rm P}(\by^k)(\bd^k)=-\mathcal{F}_{\rm P}(\by^k) 
\quad\text{ and update } \quad \by^{k+1}=\by^k+\bd^k.
\end{equation}
Notice that at each iteration of \eqref{eq:SSN_prec} one needs to evaluate
the residual function $\mathcal{F}_{\rm P}(\by^k)$,
which requires the (parallel) solution of the $N$ subproblems \eqref{eq:sub1}. 
The computational cost is therefore equivalent to one iteration of the OSM \eqref{eq:OSM_sub1}.
As an inner solver for the subproblems, which involve the (mildly) non-differentiable function $\mu$, 
a semi-smooth Newton can be employed.

We now discuss the problem of solving the Jacobian linear system in \eqref{eq:SSN_prec}. 
Let $\bd=(\bd_1,\dots,\bd_N)$, where $\bd_j=(d_{y,j},d_{p,j})$, $j=1,\dots,N$. Then a direct
calculation (omitted for brevity) shows that the
action of the operator $D\mathcal{F}_{\rm P}(\by)$ on the vector $\bd$ is given by
$D\mathcal{F}_{\rm P}(\by)(\bd) = \bd - \widetilde{\by}(\bd)$, where 
$\widetilde{\by}:=(\widetilde{\by}_1,\dots,\widetilde{\by}_N)$, and each
$\widetilde{\by}_j=(\widetilde{y}_j,\widetilde{p}_j)$  satisfies the \emph{linearized}
subdomain problems
\begin{subequations}\label{eq:lin_sub1}
\begin{align}
-\Delta \yw_j + c \yw_j + b \varphi'(y_j)\yw_j &= D\mu(p_j)(\pw_j) &&\text{ in $\Omega_j$}, \label{eq:lin_sub1_1} \\
-\Delta \pw_j + c \pw_j + b \varphi''(y_j)[p_j,\yw_j] &= \yw_j    &&\text{ in $\Omega_j$} \label{eq:lin_sub1_2} \\
\yw_j &=0, \pw_j =0  &&\text{ on $\Gamma_j^{\rm ext}$}, \label{eq:lin_sub1_3} \\
q \, \yw_j + \partial_{x} \yw_j &= q \,d_{y,j+1} + \partial_{x} d_{y,j+1} &&\text{ on $\Gamma_j$}, \label{eq:lin_sub1_4} \\
 q \, \pw_j + \partial_{x} \pw_j &= q \,d_{p,j+1} + \partial_{x} d_{p,j+1} &&\text{ on $\Gamma_j$},  \label{eq:lin_sub1_5} \\ 
 q \, \yw_j -\partial_{x} \yw_j &= q \,d_{y,j-1} - \partial_{x} d_{y,j-1} &&\text{ on $\Gamma_{j-1}$}, \label{eq:lin_sub1_6} \\
q \, \pw_j -\partial_{x} \pw_j &= q \,d_{p,j-1} - \partial_{x} d_{p,j-1} &&\text{ on $\Gamma_{j-1}$},  \label{eq:lin_sub1_7}
\end{align}
\end{subequations}
where 
\begin{equation*}
\begin{split}
D\mu(p)(\pw) = &\frac{1}{\nu}\Bigl[- \mathcal{G}_{\max}(- \beta - p) - \mathcal{G}_{\min}(\beta - p)\\
&+ \mathcal{G}_{\max}(- p - \beta-\nu \bar{u}) + \mathcal{G}_{\min}(- p + \beta + \nu \bar{u}) \Bigr] \pw,
\end{split}
\end{equation*}
with
$\mathcal{G}_{\max}(v)(x) = \begin{cases}
1 &\text{if $v(x) > 0$}, \\
0 &\text{if $v(x) \leq 0$}, \\
\end{cases}$ and
$\mathcal{G}_{\min}(v)(x) = \begin{cases}
1 &\text{if $v(x) \leq 0$}, \\
0 &\text{if $v(x) > 0$}, \\
\end{cases}$\\[10pt]
and where the boundary values for $j\in\{1,N\}$ have to be modified as in \eqref{eq:sub1}.

Note that this is the same linearized problem that must be solved repeatedly within the inner 
iterations of semi-smooth Newton, so its solution cost is only a fraction of the cost
required to calculate $\mathcal{F}_{\rm P}(\by)$.

We are now ready to state our matrix-free preconditioned semismooth Newton algorithm that corresponds
to the Newton procedure \eqref{eq:SSN_prec}.

\begin{algorithm}
\caption{Matrix-free preconditioned generalized Newton method}
\begin{algorithmic}[1]
	\begin{small}
	\REQUIRE Initial guess $\by^0$, tolerance $\epsilon$, maximum number of iterations $k_{\max}$.
	\STATE Compute $S_1(\by_2^0)$, $S_j(\by_{j-1}^0,\by_{j+1}^0)$, $j=2,\dots,N-1$, and $S_N(\by_{N-1}^0)$.
	\STATE Set $k=0$ and assemble $\mathcal{F}_{\rm P}(\by^0)$ using \eqref{eq:root_prec}.
	\WHILE{$\|\mathcal{F}_{\rm P}(\by^k)\| \geq \epsilon$ and $k \leq k_{\max}$}
		\STATE Compute $\bd^k$ by solving $D\mathcal{F}_{\rm P}(\by^k)(\bd^k)=-\mathcal{F}_{\rm P}(\by^k)$ using a matrix-free Krylov method, e.g., GMRES (together with a routine for solving \eqref{eq:lin_sub1}
		to compute the action of $D\mathcal{F}_{\rm P}(\by^k)$ on a vector $\bd$).
		\STATE Update $\by^{k+1}=\by^k+\bd^k$.
		\STATE Set $k=k+1$.
		\STATE Compute $S_1(\by_2^k)$, $S_j(\by_{j-1}^k,\by_{j+1}^k)$, $j=2,\dots,N-1$, and $S_N(\by_{N-1}^k)$.
		\STATE Assemble $\mathcal{F}_{\rm P}(\by^k)$ using \eqref{eq:root_prec}.
	\ENDWHILE
	\STATE {\bf Output}: $\by^k$.
	\end{small}
\end{algorithmic}
\label{alg:PSSN}
\end{algorithm}



\section{Numerical experiments}\label{sec:numerics}
\vspace*{-4mm}

\begin{figure}[t]
\centering
\includegraphics[scale=0.23]{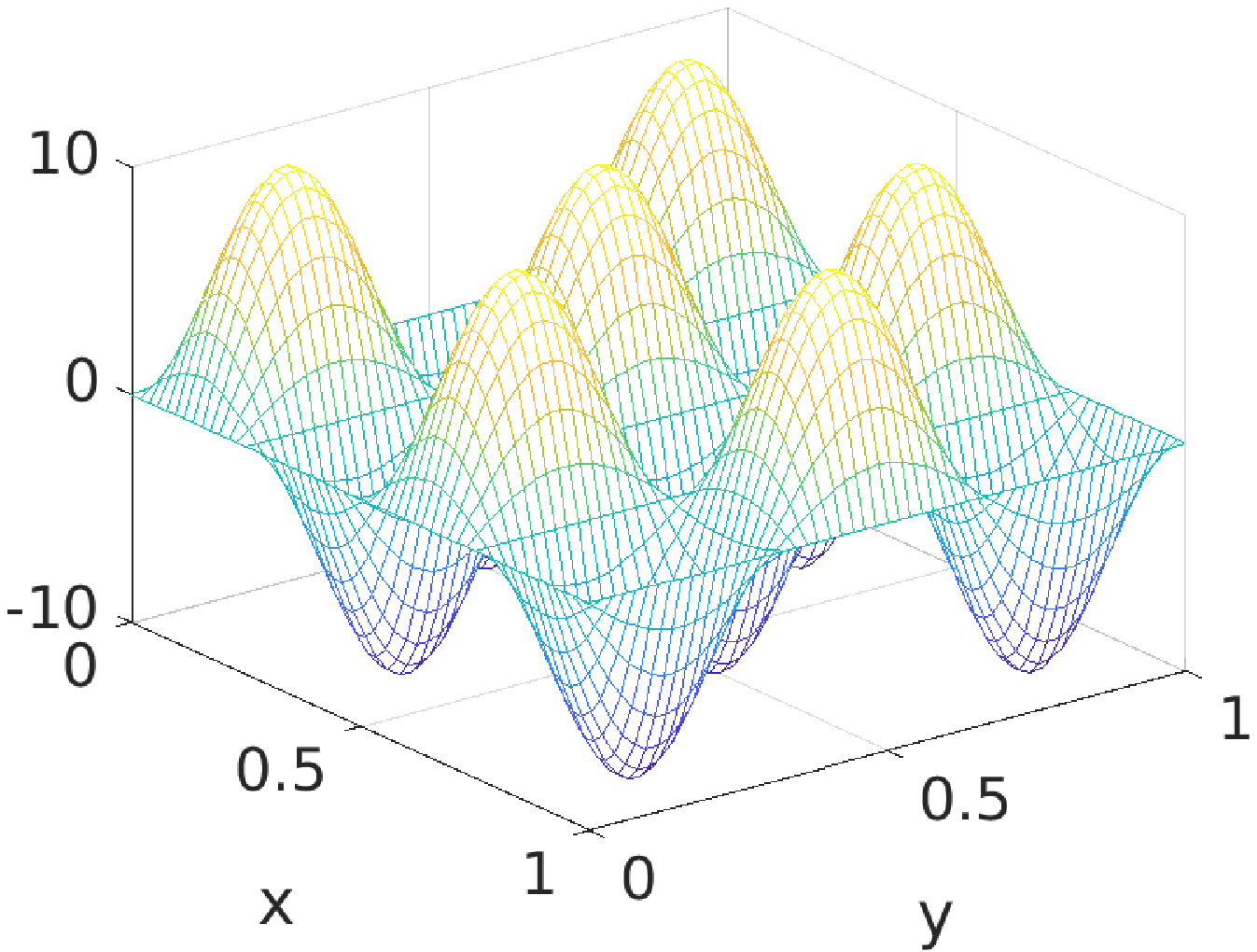}
\includegraphics[scale=0.23]{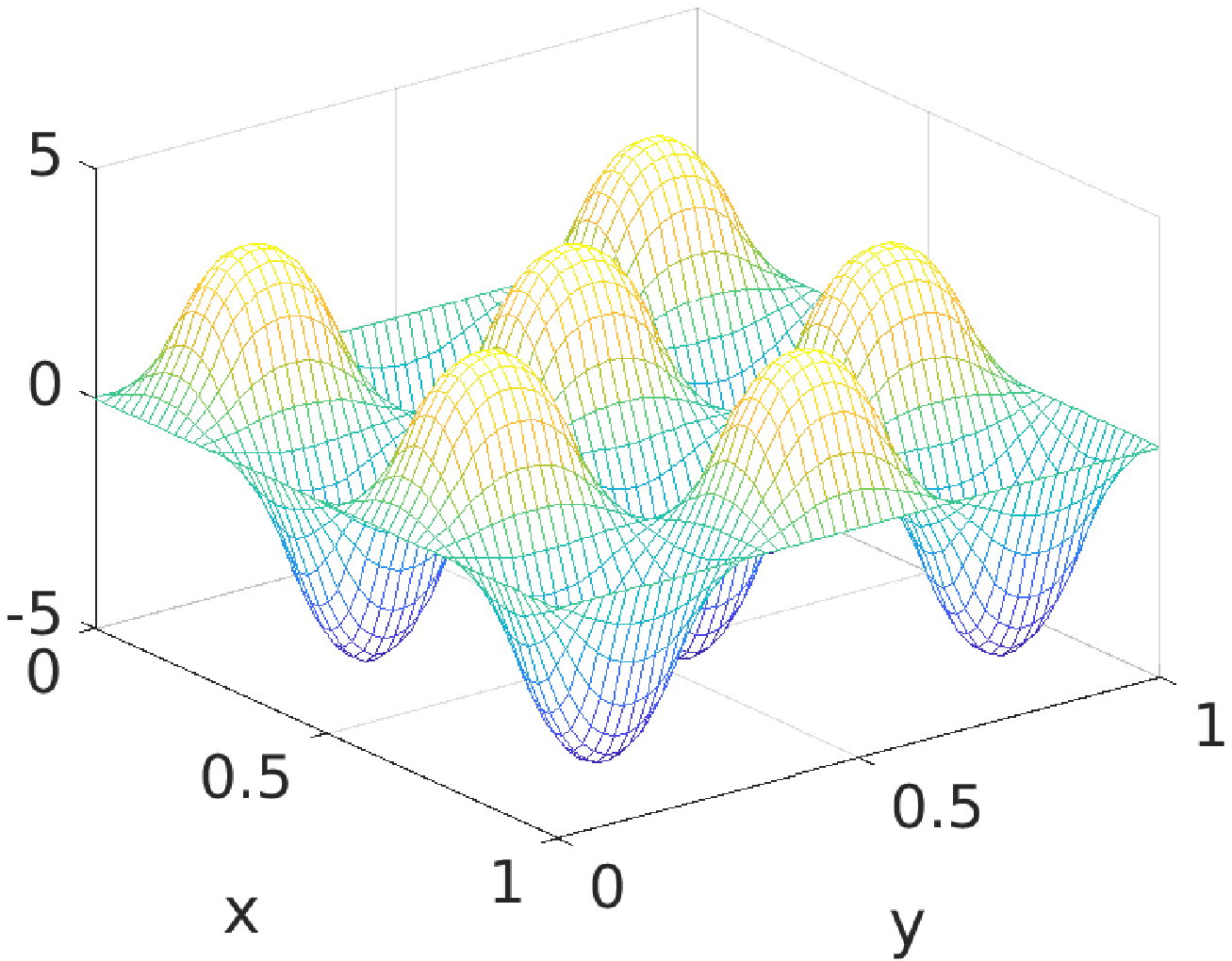}
\includegraphics[scale=0.23]{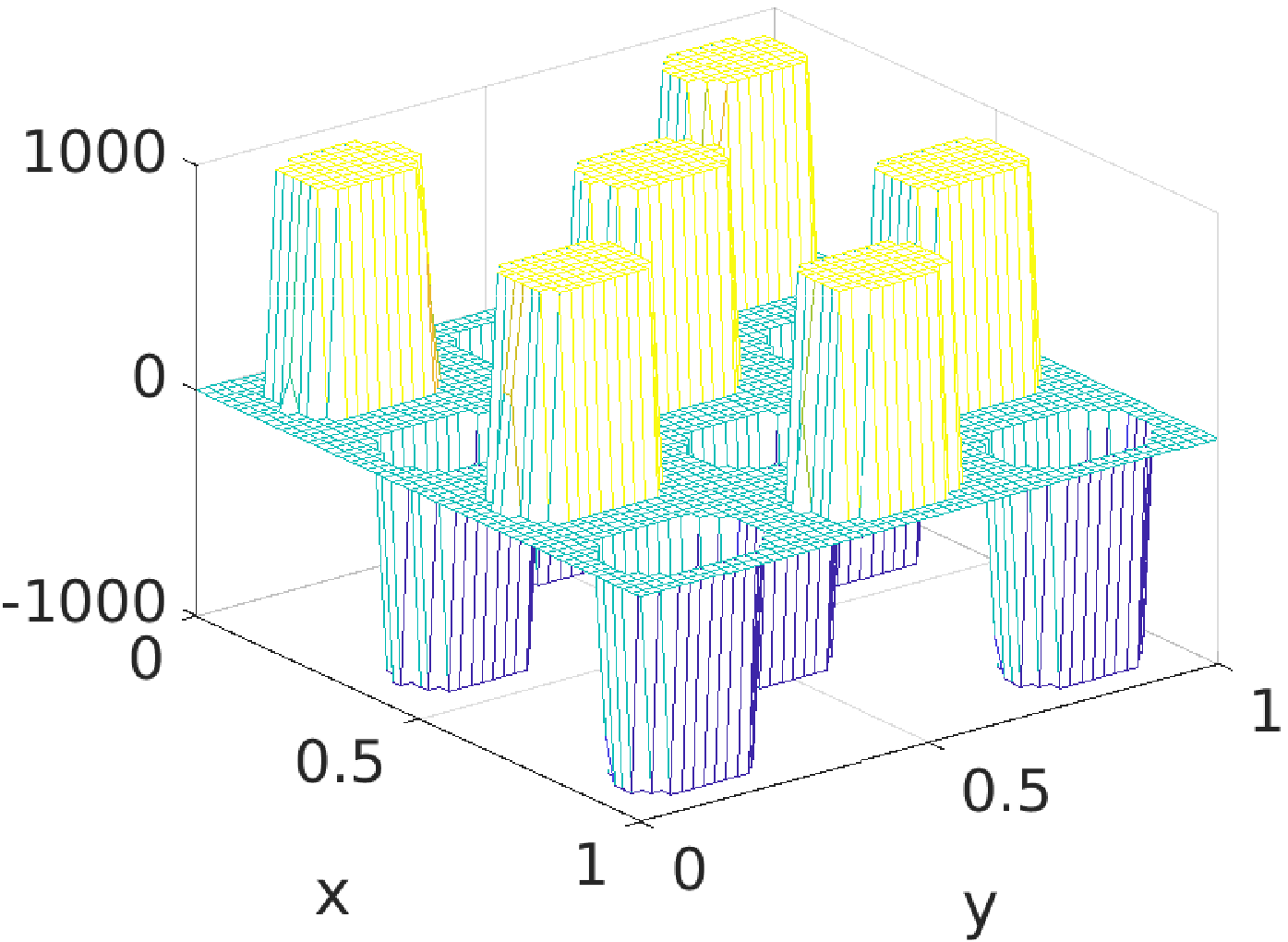}
\caption{Target $y_d$ (left), optimal state $y$ (middle), and optimal control $u$ (right) computed for $b=10$, $\nu=10^{-7}$ and $\beta = 10^{-2}$.}
\label{fig:bobo}
\end{figure}

In this section, we present results of numerical experiments.
Let us begin with a two subdomain case for 
$\Omega = (0,1)^2$, $y_d(x,y)=10 \sin(4 \pi x)\sin(3\pi y)$, 
$f=0$, $c=1$ and $\varphi(y)=y+\exp(y)$.
The domain $\Omega$ is discretized with a uniform mesh of $51$
interior points on each edge of the unit square. 
The discrete optimality system is obtained by the finite difference method.
An example of the solution computed for $b=10$, $\nu=10^{-7}$,
$\bar u=10^{3}$ and $\beta = 10^{-2}$ is shown in Fig. \ref{fig:bobo}.
Here, we can observe how the computed optimal state (middle)
has the same shape as the target $y_d$ (left).
Even though the regularization parameter $\nu$ is quite small,
the control constraints and the $L^1$-penalization prevent the control
function from making the state equal to the desired target.

To study the efficiency and the robustness of the proposed numerical framework,
we test the nonlinearly preconditioned Newton for several values of parameters $\nu$,
$\beta$, $\bar{u}$, $b$ and $q$, and compare the obtained number of iterations
with the ones performed by a (damped) semismooth Newton
applied directly to \eqref{eq:system}.
Moreover, to improve the robustness of our preconditioned Newton method,
we implemented a continuation procedure with respect to 
the regularization parameter $\nu$. This parameter is reduced over successive iterations according
to the rule 
$\nu_{k+1} = \max\{ \nu_k / 4 \, , \, \nu \}$, where 
$\nu_1=10^{-1}$ and $\nu$ is the desired final value;
see, e.g., \cite{CiaramellaNewton} for convergence results
about similar continuation procedures. 
We initialize the three methods by randomly chosen vectors. 
The number of iterations performed by both methods to reach a tolerance of
$10^{-8}$ are reported in Tab.~\ref{tab:1}, where 
the symbol $\times$ indicates non-convergence. 
From these results, 
it is clear that if the preconditioned Newton converges, then it 
 outperforms the semismooth Newton
applied directly to the full system \eqref{eq:system}. 
However, the preconditioned Newton does not always converge due to the lack of damping.
Choosing a large Robin parameter improves the robustness of the iterations,
but it is not capable of fully resolving this issue.
The continuation strategy, on the other hand, always leads to convergence 
with an iteration count comparable (for
moderate values of $\nu$) or much lower (for small values of $\nu$)
than that of the semismooth Newton method.

To better gauge the computational cost of the continuation strategy, we show the total number
of inner iterations required by `pure' preconditioned Newton versus the one with continuation in Tab.~\ref{tab:2}. The reported 
numbers are computed as $\sum_k \max_{j=1,2} {\rm it}_{j,k}$, where $k$ is the iteration count
and ${\rm it}_{j,k}$, $j=1,2$, are the number of inner iterations required by the two
subdomain solves performed at the $k$th outer iteration. (The max accounts for the fact that the
two subdomain problems are supposed to be solved in parallel.) The results show that the continuation
procedure actually \emph{reduces} the total number of inner iterations for the most part, except for
some very easy cases, such as $\beta = b = 0$, $\bar{u}=\infty$ (where the problem is in fact linear).

Finally, we remark that the GMRES iteration count is  generally 5--30 times lower for the preconditioned Newton methods than for semi-smooth
Newton applied directly to \eqref{eq:system}. This is because the Jacobians of preconditioned Newton naturally include optimized Schwarz preconditioning, whereas semismooth Newton on \eqref{eq:system}
requires additional preconditioning in order to be competitive.
%
All these numerical observation show clearly
the efficiency of the 
proposed computational framework.

\begin{table}[t]
  \centering

\parbox{\linewidth}{
\centering
\begin{tabular}{ | c | c | c | c | c | c |  c | c | c |}
  \hline
  & & & \multicolumn{3}{c|}{$\bar{u} = 10^3$} & \multicolumn{3}{c|}{$\bar{u} = \infty$}\\
  \cline{4-9}
& $q$ & $b$ & $\nu=10^{-3}$ & $\nu=10^{-5}$ & $\nu=10^{-7}$  &                         $\nu=10^{-3}$ & $\nu=10^{-5}$ & $\nu=10^{-7}$ \\ \hline                      
& 1 & 0 & 4 - 5 - 2 & 6 - 9 - 14 & $\times$ - 11 - 47 & 3 - 5 - 2 & 3 - 9 - 2 & 3 - 11 - 3\\
& 10 & 0 & 4 - 5 - 2 & 6 - 9 - 14 & 7 - 11 - 47 & 3 - 5 - 2 & 3 - 8 - 2 & 3 - 11 - 3\\                        
\parbox[b][0pt][c]{8pt}{\rotatebox[origin=c]{90}{$\beta=0$}}& 100 & 0 & 3 - 5 - 2 & 6 - 9 - 14 & $\times$ - 11 - 47 & 3 - 5 - 2 & 3 - 9 - 2 & 3 - 12 - 3\\ \cline{2-9}
 & 1 & 10 & $\times$ - 6 - 4 & $\times$ - 10 - 8 & $\times$ - 12 - 23 & 6 - 6 - 4 & $\times$ - 10 - 23 & $\times$ - 15 - $\times$\\
 & 10 & 10 & 6 - 6 - 4 & $\times$ - 10 - 8 & 9 - 12 - 23 & 6 - 6 - 4 & $\times$ - 10 - 23 & $\times$ - 14 - $\times$\\
 & 100 & 10 & 4 - 6 - 4 & 6 - 10 - 8 & $\times$ - 13 - 23 & 4 - 6 - 4 & 6 - 10 - 23 & $\times$ - 13 - $\times$\\    
 \hline
& 1 & 0 & 5 - 5 - 3 & 7 - 9 - 7 & 3 - 12 - 35 & 4 - 5 - 3 & 5 - 9 - 5 & $\times$ - 12 - 8\\
 & 10 & 0 & 5 - 5 - 3 & 5 - 9 - 7 & 2 - 11 - 35 & 4 - 5 - 3 & 4 - 9 - 5 & $\times$ - 12 - 8\\                          
\parbox[b][0pt][c]{8pt}{\rotatebox[origin=c]{90}{$\beta=10^{-2}$}}& 100 & 0 & 4 - 5 - 3 & 6 - 10 - 7 & 9 - 12 - 35 & 4 - 5 - 3 & 5 - 9 - 5 & $\times$ - 12 - 8\\  \cline{2-9}                   
& 1 & 10 & $\times$ - 6 - 4 & $\times$ - 11 - 9 & $\times$ - 12 - 37 & 6 - 6 - 4 & $\times$ - 10 - 35 & $\times$ - 13 - $\times$\\
 & 10 & 10 & 5 - 6 - 4 & $\times$ - 11 - 9 & $\times$ - 13 - 37 & 6 - 6 - 4 & $\times$ - 10 - 35 & $\times$ - 14 - $\times$\\
 & 100 & 10 & 4 - 6 - 4 & 6 - 11 - 9 & 11 - 13 - 37 & 4 - 6 - 4 & 6 - 10 - 35 & $\times$ - 13 - $\times$\\       \hline           
\end{tabular}
}

\caption{Two subdomains: outer iterations of preconditioned Newton (left value),
preconditioned Newton with continuation (middle value) and
semismooth Newton applied to the original problem (right value).
}
\label{tab:1}
\end{table}

\begin{table}[t]
\centering

\parbox{\linewidth}{
\centering
\begin{tabular}{ | c | c | c | c | c | c |  c | c | c |}
  \hline
  & & & \multicolumn{3}{c|}{$\bar{u} = 10^3$} & \multicolumn{3}{c|}{$\bar{u} = \infty$}\\
  \cline{4-9}
& $q$ & $b$ & $\nu=10^{-3}$ & $\nu=10^{-5}$ & $\nu=10^{-7}$  &   $\nu=10^{-3}$ & $\nu=10^{-5}$ & $\nu=10^{-7}$ \\ \hline                      
 & 1 & 0 & 6 - 5 & 31 - 12 & $\times$ - 18 & 2 - 5 & 3 - 8 & 3 - 11\\
 & 10 & 0 & 5 - 5 & 26 - 11 & 96 - 19 & 2 - 5 & 3 - 8 & 3 - 11\\                     
\parbox[b][0pt][c]{8pt}{\rotatebox[origin=c]{90}{$\beta=0$}}& 100 & 0 & 2 - 5 & 18 - 13 & $\times$ - 19 & 2 - 5 & 2 - 8 & 3 - 11\\ \cline{2-9}
 & 1 & 10 & $\times$ - 17 & $\times$ - 35 & $\times$ - 47 & 27 - 17 & $\times$ - 34 & $\times$ - 60\\
 & 10 & 10 & 21 - 14 & $\times$ - 31 & 103 - 43 & 21 - 14 & $\times$ - 32 & $\times$ - 53\\
 & 100 & 10 & 8 - 14 & 26 - 32 & $\times$ - 43 & 8 - 14 & 45 - 30 & $\times$ - 47\\ 
 \hline
 & 1 & 0 & 13 - 8 & 32 - 16 & 84 - 25 & 8 - 8 & 10 - 14 & $\times$ - 25\\
 & 10 & 0 & 10 - 8 & 22 - 17 & 33 - 23 & 7 - 8 & 11 - 15 & $\times$ - 24\\                        
\parbox[b][0pt][c]{8pt}{\rotatebox[origin=c]{90}{$\beta=10^{-2}$}} & 100 & 0 & 7 - 6 & 15 - 15 & 104 - 20 & 7 - 6 & 12 - 13 & $\times$ - 22\\  \cline{2-9}                   
 & 1 & 10 & $\times$ - 17 & $\times$ - 33 & $\times$ - 45 & 28 - 17 & $\times$ - 32 & $\times$ - 47\\
 & 10 & 10 & 20 - 14 & $\times$ - 33 & $\times$ - 48 & 20 - 14 & $\times$ - 30 & $\times$ - 46\\
 & 100 & 10 & 10 - 14 & 23 - 30 & 125 - 44 & 10 - 14 & 40 - 26 & $\times$ - 44\\      \hline           
\end{tabular}
}

\caption{Two subdomains: total number of inner iterations of preconditioned Newton (left value) and 
preconditioned Newton with continuation (right value).
}
\label{tab:2}
\end{table}

Let us now consider a multiple subdomain case.
In this case, the discretization mesh is refined
to have $101$ interior points on each edge of $\Omega$.
We then repeat the experiments presented above, but we fix $q=100$
and consider different numbers of subdomains, namely $N=4,8,16$. 
The results of these experiments are reported in Tab. \ref{tab:3} and \ref{tab:4},
where the number of iterations of the preconditioned Newton method (without and with
continuation) are compared to those of the semismooth Newton method applied
to \eqref{eq:system}. 
These results show that the preconditioned Newton methods are
robust against the number of subdomains, even though the size of the
subdomains decreases; see \cite{Chaouqui2018,CiaramellaGander} 
for related scalability discussions. 
Moreover, as for the two-subdomain case, the continuation procedure 
exhibits convergence in all cases.

\begin{table}[t]
\centering

\parbox{\linewidth}{
\centering
\begin{tabular}{ | c | c | c | c | c | c |  c | c | c |}
  \hline
  & & & \multicolumn{3}{c|}{$\bar{u} = 10^3$} & \multicolumn{3}{c|}{$\bar{u} = \infty$}\\
  \cline{4-9}
& $N$ & $b$ & $\nu=10^{-3}$ & $\nu=10^{-5}$ & $\nu=10^{-7}$  & $\nu=10^{-3}$ & $\nu=10^{-5}$ & $\nu=10^{-7}$ \\ \hline                      
 & 4 & 0 & 3 - 5 - 2 & 7 - 10 - 7 & $\times$ - 11 - 27 & 3 - 5 - 2 & 3 - 9 - 2 & 3 - 11 - 3\\
 & 8 & 0 & 3 - 5 - 2 & $\times$ - 10 - 7 & $\times$ - 11 - 27 & 3 - 5 - 2 & 3 - 9 - 2 & 3 - 12 - 3\\
\parbox[b][0pt][c]{8pt}{\rotatebox[origin=c]{90}{$\beta=0$}} & 16 & 0 & 3 - 5 - 2 & $\times$ - 10 - 7 & $\times$ - 11 - 27 & 3 - 5 - 2 & 3 - 9 - 2 & 3 - 12 - 3\\ \cline{2-9}
 & 4 & 10 & 4 - 6 - 4 & 7 - 11 - 7 & 10 - 14 - 21 & 4 - 6 - 4 & 6 - 10 - 10 & $\times$ - 13 - $\times$\\
 & 8 & 10 & 4 - 6 - 4 & $\times$ - 11 - 7 & $\times$ - 14 - 21 & 5 - 6 - 4 & 6 - 10 - 10 & $\times$ - 13 - $\times$\\
 & 16 & 10 & 4 - 6 - 4 & 9 - 11 - 7 & $\times$ - 14 - 21 & 4 - 6 - 4 & 6 - 10 - 10 & $\times$ - 13 - $\times$\\  
 \hline
 & 4 & 0 & 4 - 6 - 3 & 6 - 10 - 6 & 11 - 12 - 15 & 4 - 6 - 3 & 5 - 9 - 5 & 8 - 13 - 9\\
 & 8 & 0 & 4 - 6 - 3 & $\times$ - 10 - 6 & $\times$ - 12 - 15 & 4 - 6 - 3 & 6 - 9 - 5 & 8 - 13 - 9\\
\parbox[b][0pt][c]{8pt}{\rotatebox[origin=c]{90}{$\beta=10^{-2}$}}& 16 & 0 & 4 - 6 - 3 & 8 - 10 - 6 & $\times$ - 12 - 15 & 4 - 6 - 3 & 6 - 9 - 5 & 10 - 13 - 9\\  \cline{2-9}                   
 & 4 & 10 & 4 - 6 - 4 & 6 - 10 - 6 & 12 - 13 - 17 & 4 - 6 - 4 & 6 - 10 - 11 & $\times$ - 13 - $\times$\\
 & 8 & 10 & 4 - 6 - 4 & $\times$ - 11 - 6 & $\times$ - 16 - 17 & 5 - 6 - 4 & 7 - 10 - 11 & $\times$ - 13 - $\times$\\
 & 16 & 10 & 4 - 6 - 4 & 8 - 11 - 6 & $\times$ - 18 - 17 & 4 - 6 - 4 & 8 - 10 - 11 & $\times$ - 14 - $\times$\\ \hline           
\end{tabular}
}

\caption{Multiple subdomains: outer iterations of preconditioned Newton (left value),
preconditioned Newton with continuation (middle value)
and semismooth Newton applied to the original system (right value).
}
\label{tab:3}
\end{table}

\begin{table}[t]
\centering
\parbox{\linewidth}{
\centering
\begin{tabular}{ | c | c | c | c | c | c |  c | c | c |}
  \hline
  & & & \multicolumn{3}{c|}{$\bar{u} = 10^3$} & \multicolumn{3}{c|}{$\bar{u} = \infty$}\\
  \cline{4-9}
& $N$ & $b$ & $\nu=10^{-3}$ & $\nu=10^{-5}$ & $\nu=10^{-7}$  &   $\nu=10^{-3}$ & $\nu=10^{-5}$ & $\nu=10^{-7}$ \\ \hline                      
 & 4 & 0 & 2 - 5 & 21 - 15 & $\times$ - 18 & 2 - 5 & 2 - 8 & 3 - 11\\
 & 8 & 0 & 2 - 5 & $\times$ - 12 & $\times$ - 15 & 2 - 5 & 2 - 8 & 3 - 11\\                   
\parbox[b][0pt][c]{8pt}{\rotatebox[origin=c]{90}{$\beta=0$}} & 16 & 0 & 4 - 5 & $\times$ - 13 & $\times$ - 15 & 2 - 5 & 2 - 8 & 3 - 11\\ \cline{2-9}
 & 4 & 10 & 8 - 12 & 32 - 30 & 75 - 47 & 9 - 12 & 31 - 28 & $\times$ - 45\\
 & 8 & 10 & 8 - 11 & $\times$ - 27 & $\times$ - 44 & 8 - 11 & 28 - 26 & $\times$ - 42\\
 & 16 & 10 & 7 - 11 & 29 - 27 & $\times$ - 36 & 7 - 11 & 27 - 24 & $\times$ - 39\\
 \hline
 & 4 & 0 & 7 - 8 & 17 - 19 & 61 - 22 & 7 - 8 & 11 - 15 & 44 - 29\\
 & 8 & 0 & 7 - 7 & $\times$ - 15 & $\times$ - 19 & 7 - 7 & 13 - 14 & 35 - 26\\                      
\parbox[b][0pt][c]{8pt}{\rotatebox[origin=c]{90}{$\beta=10^{-2}$}}  & 16 & 0 & 6 - 6 & 14 - 15 & $\times$ - 18 & 5 - 6 & 9 - 12 & 32 - 24\\  \cline{2-9}                   
 & 4 & 10 & 10 - 13 & 23 - 27 & 106 - 44 & 10 - 13 & 31 - 27 & $\times$ - 43\\
 & 8 & 10 & 10 - 11 & $\times$ - 27 & $\times$ - 46 & 9 - 11 & 30 - 26 & $\times$ - 42\\
 & 16 & 10 & 8 - 11 & 24 - 27 & $\times$ - 44 & 8 - 11 & 30 - 23 & $\times$ - 38\\    \hline           
\end{tabular}
}
\caption{Multiple subdomains: total number of inner iterations of preconditioned Newton (left value) and
preconditioned Newton with continuation (right value).
}
\label{tab:4}
\end{table}

\section{Further discussion and conclusion}\label{sec:conclusion}
\vspace*{-4mm}
This short manuscript represents a proof of concept for 
using domain decom\-position-based nonlinear preconditioning 
to efficiently solve nonlinear, nonsmooth optimal control problems
governed by elliptic equations. However, several theoretical and numerical issues
must be addressed as part of 
a complete development of these techniques.
From a theoretical point of view, to establish concrete convergence results based on 
classical semismooth Newton theory, it is crucial to study the (semismoothness)
properties of the subdomain solution maps $\mathcal{S}_j$, which are implicit function
of semi-smooth maps. 
Another crucial point is the proof of well-posedness of the
(preconditioned) Newton linear system.
From a 
domain decomposition perspective, more general decompositions
(including cross points) must be considered. 
Finally, a detailed analysis of the scalability of the GMRES iterations is necessary.

	\bibliographystyle{plain}
	\vspace*{-4mm}
	\bibliography{references}
\end{document}